\documentstyle[12pt]{article}
\makeatletter
\renewcommand\thesection{\@Roman\c@section}
\renewcommand\thesubsection{\thesection.\@arabic\c@subsection}
\makeatother

\begin{document}
\begin{titlepage}
\begin{flushright}
math.QA/9911058
\end{flushright}
\vskip.3in

\begin{center}
{\Large\bf $U_q[\widehat{sl(2|1)}]$ Vertex Operators, Screen Currents and
Correlation Functions at Arbitrary Level}
\vskip.3in
{\large Yao-Zhong Zhang} and {\large Mark D. Gould}
\vskip.2in
{\em Department of Mathematics, University of Queensland, Brisbane,
     Qld 4072, Australia

Email: yzz@maths.uq.edu.au}
\end{center}

\vskip 2cm
\begin{center}
{\bf Abstract}
\end{center}

Bosonized $q$-vertex operators related to the 4-dimensional evaluation 
modules of the quantum affine superalgebra 
$U_q[\widehat{sl(2|1)}]$ are constructed for arbitrary
level $k=\alpha$, where $\alpha\neq 0, -1$ is a complex parameter
appearing in the 4-dimensional evaluation representations. They are
intertwiners among the level-$\alpha$ highest weight Fock-Wakimoto 
modules. Screen currents which commute
with the action of $U_q[\widehat{sl(2|1)}]$ up to total differences
are presented. Integral formulae for $N$-point functions of type I
and type II $q$-vertex operators are proposed.


\end{titlepage}


\def\a{\alpha}
\def\b{\beta}
\def\d{\delta}
\def\e{\epsilon}
\def\ve{\varepsilon}
\def\g{\gamma}
\def\k{\kappa}
\def\l{\lambda}
\def\o{\omega}
\def\t{\theta}
\def\s{\sigma}
\def\D{\Delta}
\def\L{\Lambda}

\def\F{{\cal F}}
\def\G{{sl(2|1)}}
\def\hG{{\widehat{sl(2|1)}}}
\def\R{{\cal R}}
\def\hR{{\hat{\cal R}}}
\def\C{{\bf C}}
\def\P{{\bf P}}
\def\Z2{{{\bf Z}_2}}
\def\T{{\cal T}}
\def\H{{\cal H}}
\def\trho{{\tilde{\rho}}}
\def\tphi{{\tilde{\phi}}}
\def\tT{{\tilde{\cal T}}}
\def\Uq{U_q[\widehat{sl(2|1)}]}


\def\beq{\begin{equation}}
\def\eeq{\end{equation}}
\def\bea{\begin{eqnarray}}
\def\eea{\end{eqnarray}}
\def\ba{\begin{array}}
\def\ea{\end{array}}
\def\no{\nonumber}
\def\lt{\left}
\def\rt{\right}
\newcommand{\bq}{\begin{quote}}
\newcommand{\eq}{\end{quote}}

\newtheorem{Theorem}{Theorem}
\newtheorem{Definition}{Definition}
\newtheorem{Conjecture}{Conjecture}
\newtheorem{Proposition}{Proposition}
\newtheorem{Lemma}{Lemma}
\newtheorem{Corollary}{Corollary}
\newcommand{\proof}[1]{{\bf Proof. }
        #1\begin{flushright}$\Box$\end{flushright}}

\newcommand{\sect}[1]{\setcounter{equation}{0}\section{#1}}
\renewcommand{\theequation}{\thesection.\arabic{equation}}

\sect{Introduction\label{intro}}

The notion of $q$-vertex operators as certain intertwiners of highest
weight modules of quantum affine algebras was introduced by Frenkel
and Reshetikhin \cite{Fre92} in their work on the $q$-deformation of
the Wess-Zumino-Novikov-Witten model. These $q$-vertex operators give
rise to $q$-analogues of the primary fields in conformal field theory.

Similar to the classical case, $q$-vertex operators are characterized
by the intertwining property defined from the relevant quantum affine algebras.
However it is non-trivial to obtain explicit expressions of them.
A powerful tool for constructing such explicit formulae
is the bosonization technique 
\cite{Fre88,Ber89,Awa94}, initiated by Wakimoto \cite{Wak86} in the
theory of affine Lie algebras. This method enables one
in principle to determine $q$-vertex operators in terms of certain
free bosonic fields. So far, level-one bosonized
$q$-vertex operators 
have been constructed for most quantum affine algebras 
\cite{Koy94,Jin95,Jin99} and the type I quantum affine superalgebras
$U_q[\widehat{sl(M|N)}],~M\neq N$ 
\cite{Kim97} and $U_q[\widehat{gl(N|N)}]$ \cite{Zha98,Yan99b}. 
In the case of arbitrary level, bosonized formulae have been known
only for the type I $q$-vertex operators of $U_q[\widehat{sl(2)}]$
\cite{Kat93,Mat94,Bou94,Kon94} and
$U_q[\widehat{sl(N)}]$ \cite{Awa94}.

One of the central issues in conformal field theory and massive integrable
models is the computation of correlation functions,
which are  matrix elements of certain products of vertex
operators. The explicit bosonized expressions of vertex operators
play an essential role. They enable one to compute
correlators exactly in the form of integral representations. This
was demonstrated by the Kyoto group and collaborators in their
ground-breaking work on the diagonalization of the $XXZ$ spin
chain \cite{Dav93,Jim94}.
In \cite{Koy94,Yan99,Hou99}, certain correlation functions of other
quantum affine (super)algebras at level one
were computed via the bosonization procedure, generalizing the
work of the Kyoto group and collaborators.

The case of arbitrary level is more complicated. Due to the existence
of nontrivial backgound charges, the naive solutions to the intertwining
relations in terms of free bosonic fields do not give rise to proper
bosonizations of the $q$-vertex operators, which ensure the nonvanishing
of correlation functions.
As in conformal field theory, q-screen currents which balance the 
background charges are generally needed.
$q$-screen currents are dimension 1 
operators which (anti-)commute with the relevant quantum algebra 
generators up to total differences.
Bosonized $q$-screen currents have been obtained for $U_q[\widehat{sl(N)}]$
\cite{Kat93,Mat94,Bou94,Kon94,Awa94} and been applied to
compute the correlation functions of the type I
$U_q[\widehat{sl(2)}]$ vertex operators \cite{Kat93,Mat94,Bou94,Kon94}.

In this paper, by using the free field realization of $\Uq$ at arbitrary level
$k\neq 0,-1$ \cite{Awa97} we investigate the bosonization of
$q$-vertex operators related to the 4-dimensional evaluation modules of
$\Uq$. It is worth mentioning that our
4-dimensional representation contains an extra complex parameter
$\a\neq 0,-1$. For arbitrary level 
$k=\a$, the $q$-vertex operators are mappings of 
certain highest weight Fock-Wakimoto modules in a bosonic Fock
space. Screen currents which (anti-)commute with the action of $\Uq$
are obtained and bosonized $q$-vertex operators dressed with the screen
charges are proposed. This provides a natural way to write down
an integral representation for correlation functions of the
bosonized $q$-vertex operators. 

The results obtained in this paper will be useful in analysing the
supersymmetric integrable model introduced in \cite{Bra95}. This is a
quantum spin chain model arising from the R-matrix for the 4-dimensional
$\Uq$ evaluation representation and can be interpreted as a model
describing strongly correlated electrons.

\sect{Prelimilaries}

\subsection{Quantum affine superalgebra $U_q[\widehat{sl(2|1)}]$}

The simple roots of the affine superalgebra $\widehat{sl(2|1)}$ 
\cite{Kac78} are
\bea
\a_0=\d-\ve_1+\d_1,~~~~\a_1-\ve_1-\ve_2,~~~~\a_2=\ve_2-\d_1,\no
\eea
where $\d$ is the null root and $\{\ve_1,\ve_2,\d_1\}$ are
orthonormal basis satisfying
\bea
&&(\d,\d)=(\d,\ve_i)=(\d,\d_1)=(\d_1,\ve_i)=0,~~~i=1,2,\no\\
&&(\ve_i,\ve_j)=\d_{ij},~~~~(\d_1,\d_1)=-1.\no
\eea
The fundamental weights are
\bea
\L_0,~~~~\L_1=\L_0-\ve_2+\d_1,~~~~\L_2=\L_0-\ve_1-\ve_2+2\d_1,\no
\eea
where $\L_0$ is the affine weight obeying
$(\L_0,\L_0)=(\L_0,\ve_i)=0,~i=1,2$ and $(\L_0,\d)=1$. 
The symmetric Cartan matrix $(a_{ij})$ 
of the affine Lie superalgebra $\widehat{sl(2|1)}$ has elements
$a_{ij}=(\a_i,\a_j),~i,j=0,1,2$. Explicitly,
\begin{eqnarray*}
(a_{ij})=\left(
\begin{array}{ccc}
0&-1&1\\
-1&2&-1\\
1&-1&0
\end{array}\right).
\end{eqnarray*}

Quantum affine superalgebra 
$U_q[\widehat{sl(2|1)}]$ is a $q$-analogue of the universal 
enveloping algebra of $\widehat{sl(2|1)}$ generated by the Chevalley 
generators $\{e_i,f_i,q^{h_i},d | i=0,1,2\}$, where $d$ is the
usual derivation operator. The $\Z2$-grading
of the generators are
$[e_0]=[f_0]=[e_{2}]=[f_{2}]=1$ and zero otherwise. 
The defining relations are
\bea
& & h_ih_j =h_jh_i,\ \ \ h_id=dh_i, \ \ \ [d,e_i]=\delta_{i,0}e_i,\ \ \
[d,f_i]=-\delta_{i,0}f_i,\no\\ 
& &q^{h_i}e_jq^{-h_i}=q^{a_{ij}}e_j,\ \ \ 
 q^{h_i}f_jq^{-h_i}=q^{-a_{ij}}f_j ,\ \ \ 
 [e_i,f_j] =\delta_{ij}{ q^{h_i}-q^{-h_i} \over q-q^{-1}},\no\\
& &[e_i, e_j]=[f_i, f_j]=0,~~~{\rm for}~a_{ij}=0,\no\\
& &[e_1, [e_1, e_l]_{q^{-1}}]_q=0,\ \ \
   [f_1, [f_1, f_l]_{q^{-1}}]_q=0,\ \ l=0,2.
\eea
Here and throughout, $[X,Y]_\xi=XY-(-1)^{[X][Y]}\xi YX$ and $[X,Y]=
[X,Y]_1$.
 
$U_q[\widehat{sl(2|1)}]$ is a quasi-triangular Hopf superalgebra
endowed with the ${\bf Z}_2$-graded Hopf algebra structure:
\bea
&&\Delta(h_i)=h_i\otimes 1+1\otimes h_i,~~~~
  \D(d)=d\otimes 1+1\otimes d,\no\\
&& \Delta(e_i)=e_i\otimes 1+q^{h_i}\otimes e_i,~~~~
  \Delta(f_i)=f_i\otimes q^{-h_i}+1\otimes f_i ,\no\\
&&\epsilon(h_i)=\e(d)=\epsilon(e_i)=\epsilon(f_i)=0,\no\\
&&S(e_i)=-q^{-h_i} e_i,~~~~ S(f_i)=-f_i q^{h_i}, ~~~~
   S(h_i)=-h_i,~~~~ S(d)=-d.
\eea 
Note the antipode $S$ is a $\Z2$-graded algebra anti-automorphism. Namely for
homogeneous elemets $a,b \in \Uq$,
$S(ab)=(-1)^{[a][b]}S(b)S(a)$.
The multiplication rule for the tensor product is ${\bf Z}_2$
graded and is defined for homogeneous elements $a,b,a',b' \in \Uq$ by 
$(a\otimes b)(a'\otimes b')=(-1)^{[b][a']}(aa'\otimes bb')$,  
which extends to inhomgeneous elements through linearity.

$\Uq$ can also be realized by 
the Drinfeld generators \cite{Dri88} $\{X^{\pm,i}_{m}$,
$h^{i}_n$, $q^{h_0^i}, c, d | i=1,2,$  $m \in 
{\bf Z}, n \in {\bf Z}_{\ne 0}\}$. 
The ${\bf Z}_2$-grading of the Drinfeld generators 
are $[X^{\pm,2}_m]=1$ ($m\in {\bf Z}$) and zero otherwise. The
relations read \cite{Yam96,Zha97}
\bea
& &c: \  {\rm central ~ element},\no\\
& &[h_0^i,h^j_m]=0,\ \ \ [d,h_0^i]=0,\ \ \ [d,h^{j}_m]=mh^j_m,\no\\
& &[h^i_m,h^j_n] =\delta_{m+n,0} 
\frac{[a_{ij}m]_q[nc]_q}{m},\no \\
& &q^{h_0^i}X^{\pm,j}_mq^{-h^i_0} =q^{\pm a_{ij}}X^{\pm,j}_m,\ \ \
[d,X^{\pm,j}_m]=mX^{\pm,j}_m,\no \\ 
& &[h^i_m,X^{\pm,j}_n]=\pm {[a_{ij}m]_q \over m}
q^{\pm |m|c/2}X^{\pm,j}_{n+m},\no\\ 
& &[X^{+,i}_m,X^{-,j}_n]=\frac{\delta_{i,j}}{ q-q^{-1}}
(q^{(m-n)c/2}\psi^{+,j}_{m+n} -q^{-(m-n)c/2 }
\psi^{-,j}_{m+n}),\no\\ 
& &[X^{\pm ,2}_m,X^{\pm ,2}_n]=0,\no \\
& &[X^{\pm,i}_{m+1},X^{\pm,j}_{n}]_{q^{\pm a_{ij}}} 
+[X^{\pm,j}_{n+1},X^{\pm,i}_{m}]_{q^{\pm a_{ij}}}=0, 
\ \ {\rm for} \ \ a_{ij}\neq 0,\no\\
& &[X^{\pm,1}_{n_1}, [X^{\pm,1}_{n_2}, X^{\pm, 2}_m]_{q^{-1}}]_q
   +(n_1\leftrightarrow n_2)=0,
\eea 
where $[m]_q=\frac{q^m-q^{-m}}{q-q^{-1}}$ and $\psi^{\pm,i}_n$ are
defined by
\bea
\sum_{n\in {\bf Z}}\psi^{\pm,i}_n z^{-n}=q^{\pm h^i_0}\exp
  \lt(\pm(q-q^{-1})\sum_{n>0}h^i_{\pm n} z^{\mp n}\rt).\no
\eea
The Chevalley generators are related to the Drinfeld generators
by the formulae:
\begin{eqnarray}
& &h_i= h_0^i,\ \ \ e_i = X^{+,i}_0,\ \ \ h_0=c-h^1_0-h^2_0,
   \ \ \ f_i = X^{-,i}_0, \ \ i=1,2,\no\\ 
& &e_0= -[X_0^{-,2},X^{-,1}_1]_{q^{-1}}q^{-h^1_0-h^2_0},\ \ \
   f_0=q^{h^1_0+h^2_0}[X^{+,1}_{-1},X^{+,2}_0]_q.
\end{eqnarray}

\subsection{Bosonization of $\Uq$ at arbitrary level $k$}

In this subsection we briefly recall the free boson realization of
$\Uq$ at arbitrary level $k$ \cite{Awa97}.
Let us introduce the bosonic $q$-oscillators
$\{a^1_n,a^2_n,b^{ij}_n,c_n$, $Q_{a^1},Q_{a^2}, Q_{b^{ij}},Q_{c}$ 
$|n \in {\bf Z}, 1\leq i<j\leq 3 \}$ which satisfy the commutation relations
\bea
& &[a^i_m,a^j_n]=\delta_{m+n,0}{[a_{ij}m]_q[(k+1)m]_q \over m},\ \ \  
[a^i_0,Q_{a^j}] = (k+1)a_{i,j},\no\\
& &[b^{ij}_m,b^{i'j'}_n]=(-1)^{\d_{j2}}\d^{ii'}\d^{jj'}
   \delta_{m+n,0}{[m]^2_q \over m},\ \ \  
[b^{ij}_0,Q_{b^{i'j'}}] = (-1)^{\d_{j2}}\d^{ii'}\d^{jj'}, \no \\
& &[c_m,c_n]=\delta_{m+n,0}{[m]^2_q \over m},\ \ \ \  
[c_0,Q_{c}] = 1. 
\eea
The remaining commutators vanish. Here and throughout
$k\neq 0, -1$ is a complex
parameter. For any pair $(a_n, Q_a)$, we define
\bea
&&a(z;\k)=-\sum_{n \ne 0}{a_n \over 
  [n]_q}q^{-\k|n|}z^{-n} +Q_{a}+a_0\ln z,\no\\
&&a_{\pm}(z)=\pm(q-q^{-1})\sum_{n>0} a_{\pm n} z^{\mp n}\pm a_0\ln q.
\eea
We have 

\begin{Theorem} \cite{Awa97}: 
Define the fields $X^{\pm,i}(z)$ by
\bea
X^{\pm,i}(z)=\sum_{n\in{\bf Z}} X^{\pm,i}_n
    z^{-n-1}.\no
\eea
Then at arbitrary level $k\neq 0, -1$, 
$\Uq$ is realized by the free boson fields as follows
\bea
c&=&k,~~~~h^1_0={a^1_0+2b^{12}_0+b^{13}_0-b^{23}_0},~~~~
  h^2_0={a^2_0-b^{12}_0-b^{13}_0},\no\\
h^1_m&=& a^1_mq^{-\frac{|m|}{2}}+b^{12}_mq^{-(\frac{k}{2}+1)|m|}
  (q^{|m|}+q^{-|m|})
  +b^{13}_mq^{-(\frac{k}{2}+2)|m|}-b^{23}_mq^{-(\frac{k}{2}+1)|m|},\no\\
h^2_m& =&a^2_mq^{-\frac{|m|}{2}}-b^{12}_mq^{-(\frac{k}{2}+1)|m|}-b^{13}_m
    q^{-(\frac{k}{2}+1)|m|},\no\\ 
X^{+,1}(z) &=& -\frac{1}{(q-q^{-1})z}:e^{-b^{12}(z;-1)}(
   e^{-c(qz;0)}-e^{-c(q^{-1}z;0}):e^{\sqrt{-1}\pi (c_0+b^{12}_0)},\no\\
X^{+,2}(z)& =& -:e^{-b^{12}_+(qz)-b^{13}_+(qz)+b^{23}(qz;0)}:
   e^{\sqrt{-1}\pi (c_0+b^{12}_0+b^{13}_0+b^{23}_0)}
   +:e^{b^{12}(z;0)+b^{13}(z;0)+c(z;0)}:,\no\\ 
X^{-,1}(z)&=& \frac{1}{(q-q^{-1})z}
  :\lt(e^{a^1_+(q^{\frac{k+1}{2}}z)+b^{12}(q^{k+2}z;1)
  +b^{13}_+(q^{k+2}z)-b^{23}_+(q^{k+1}z)+c(q^{k+1}z;0)}\rt.\no\\
& &\lt.-e^{a^1_-(q^{-\frac{k+1}{2}}z)+b^{12}(q^{-k-2}z;1)
  +b^{13}_-(q^{-k-2}z)-b^{23}_-(q^{-k-1}z)+c(q^{-k-1}z;0)}\rt):
  e^{-\sqrt{-1}\pi (c_0+b^{12}_0)}\no\\
& &+q^{k+1} :e^{a^1_+(\frac{k+1}{2}z)-b^{13}(q^{k+1}z;0)
  +b^{23}(q^{k+1}z;-1)}: e^{\sqrt{-1}\pi (b^{13}_0+b^{23}_0)},\no\\
X^{-,2}(z)&=&\frac{1}{(q-q^{-1})z}\no\\
& &\lt( q:(e^{a^2_+(q^{\frac{k+1}{2}}z)-b^{23}(q^{k+1}z;0)}
  -e^{a^2_-(q^{-\frac{k+1}{2}}z)-b^{23}(q^{-k-1}z;0)}):
  e^{-\sqrt{-1}\pi (c_0+b^{12}_0+b^{13}_0+b^{23}_0)}\rt.\no\\
& & \lt.-:e^{a^2_-(q^{-\frac{k+1}{2}}z)-b^{12}(q^{-k-1}z;1)
   -b^{13}(q^{-k-1}z;1)}
   (e^{-c(q^{-k}z;0)}-e^{-c(q^{-k-2}z;0)}):\rt).\no\\
\eea
\end{Theorem}

\sect{Level-zero representations}

We discuss level-zero representations of $\Uq$, which are needed in
next section for the investigation of $q$-vertex operators.

Let $V_\a$ is the one parameter family of the 
4-dimensional typical irreducible representation of
$U_q[sl(2|1)]$. Here and throughout,
$\a\neq 0, -1$ is a complex parameter. We choose the
basis vectors $\{v_1,v_2,v_3,v_4\}$ of $V_\a$ and assign them the $\Z2$
gradings $[v_1]=[v_4]=0,~ [v_2]=[v_3]=1$.
Let $e_{ij}$ be the $4\times 4$ matrices satisfying 
$(e_{ij})_{kl}=\d_{ik}\d_{jl}$.
In the homogeneous gradation, the
evaluation representation $V_{\a,z}$ of $\Uq$ is given by
\bea
&& e_1=e_{23},~~~~f_1=e_{32},~~~~
 h_1={e_{22}-e_{33}},\no\\
&&e_2=\sqrt{[\a]_q}e_{12}+\sqrt{[\a+1]_q}e_{34},\no\\
&&f_2=\sqrt{[\a]_q}e_{21}+\sqrt{[\a+1]_q}e_{43},\no\\
&&h_2={\a(e_{11}+e_{22})+(\a+1)(e_{33}+e_{44})},\no\\
&&e_0=-z(-\sqrt{[\a]_q}e_{31}+\sqrt{[\a+1]_q}e_{42}),\no\\
&&f_0=z^{-1}(-\sqrt{[\a]_q}e_{13}+\sqrt{[\a+1]_q}e_{24}),\no\\
&&h_0={-\a(e_{11}+e_{33})-(\a+1)(e_{22}+e_{44})}.
\eea

We define the dual module $V_{\a,z}^{*S}$ of $V_{\a,z}$ by 
$\pi_{V_{\a,z}^{*S}}(a)=\pi_{V_{\a,z}}(S(a))^{st}$, $\forall a\in \Uq$, where 
$st$ is the supertransposition operation. On $V^{*S}_{\a,z}$, the Chevalley
generators are represented by
\bea
&& e_1=-q^{-1}e_{32},~~~~f_1=-qe_{23},~~~~
  h_1={-e_{22}+e_{33}},\no\\
&&e_2=q^{-\a}\sqrt{[\a]_q}e_{21}-q^{-\a-1}\sqrt{[\a+1]_q}e_{43},\no\\
&&f_2=-q^\a\sqrt{[\a]_q}e_{12}+q^{\a+1}\sqrt{[\a+1]_q}e_{34},\no\\
&&h_2={-\a(e_{11}+e_{22})-(\a+1)(e_{33}+e_{44})},\no\\
&&e_0=-z(q^\a\sqrt{[\a]_q}e_{13}+q^{\a+1}\sqrt{[\a+1]_q}e_{24}),\no\\
&&f_0=-z^{-1}(q^{-\a}\sqrt{[\a]_q}e_{31}+q^{-\a-1}
   \sqrt{[\a+1]_q}e_{42}),\no\\
&&h_0={\a(e_{11}+e_{33})+(\a+1)(e_{22}+e_{44})}.
\eea

We state
\begin{Proposition}:
The Drinfeld generators are represented on $V_{\a,z}$ by
\bea
&&h^1_0=e_{22}-e_{33},~~~~h^2_0=\a(e_{11}+e_{22})+(\a+1)(e_{33}+e_{44}),\no\\
&&X^{+,1}_m=(zq^{\a+1})^m e_{23},~~~~
  X^{-,1}_m=(zq^{\a+1})^m e_{32},\no\\
&&X^{+,2}_m=(zq^{\a+1})^m(q^{-m}\sqrt{[\a]_q}e_{12}
  +q^m\sqrt{[\a+1]}e_{34}),\no\\
&&X^{-,2}_m=(zq^{\a+1})^m(q^{-m}\sqrt{[\a]_q}e_{21}
  +q^m\sqrt{[\a+1]_q}e_{43}),\no\\
&&h^1_m=(zq^{\a+1})^m\frac{[m]_q}{m}(q^{-m}e_{22}-q^me_{33}),\no\\
&&h^2_m=\frac{z^m}{m}\lt([\a m]_q(e_{11}+e_{22})+q^m[(\a+1)m]_q
  (e_{33}+e_{44})\rt),
\eea
and on $V^{*S}_{\a,z}$ by
\bea
&&h^1_0=-e_{22}+e_{33},~~~~h^2_0=-\a(e_{11}+e_{22})-(\a+1)(e_{33}+e_{44}),\no\\
&&X^{+,1}_m=-z^mq^{-m\a-m-1} e_{32},~~~~
  X^{-,1}_m=-z^mq^{-m\a-m+1} e_{23},\no\\
&&X^{+,2}_m=z^mq^{-(1+m)\a}(\sqrt{[\a]_q}e_{21}
  -q^{-2m-1}\sqrt{[\a+1]_q}e_{43}),\no\\
&&X^{-,2}_m=z^mq^{(1-m)\a}(-\sqrt{[\a]_q}e_{12}
  +q^{-2m+1}\sqrt{[\a+1]_q}e_{34}),\no\\
&&h^1_m=-(zq^{-\a-1})^m\frac{[m]_q}{m}(q^me_{22}-q^{-m}e_{33}),\no\\
&&h^2_m=-\frac{z^m}{m}\lt([\a m]_q(e_{11}+e_{22})+q^{-m}[(\a+1)m]_q
  (e_{33}+e_{44})\rt).
\eea
\end{Proposition}

\sect{Vertex operators at arbitrary level $k=\a$}

Let $V(\lambda)$ be a level-$k$ highest weight $\Uq$-module with 
highest weight $\lambda$ and highest weight vector $|\l>$. 
Consider the following intertwiners of 
$\Uq$-modules:
\bea
& &\Phi_{\lambda}^{\mu V}(z) :
 V(\lambda) \longrightarrow V(\mu)\otimes V_{\a,z} ,\ \ \ \ 
\Phi_{\lambda}^{\mu V^{*}}(z) :
 V(\lambda) \longrightarrow V(\mu)\otimes V_{\a,z}^{*S} ,\no\\
& &\Psi_{\lambda}^{V \mu}(z) :
 V(\lambda) \longrightarrow V_{\a,z}\otimes V(\mu),\ \ \ \
\Psi_{\lambda}^{V^* \mu}(z) :
 V(\lambda) \longrightarrow V_{\a,z}^{*S}\otimes V(\mu).\label{intertwiners}
\eea
They are intertwiners in the sense that for any $x\in \Uq$, 
\begin{eqnarray}
\Theta(z)\cdot x=\Delta(x)\cdot \Theta(z),\ \ \ \Theta(z)=
\Phi(z),\Phi^{*}(z),\Psi(z),\Psi^{*}(z).\label{intertwining}
\end{eqnarray}
The intertwiners are even operators, that is their grading
is $[\Theta(z)]=0$. $\Phi(z)$ 
($\Phi^{*}(z)$) is called type I (dual) vertex operator and $\Psi(z)$ 
($\Psi^{*}(z)$) type II (dual) vertex operator.

Expand these vertex operators in terms of their components 
\begin{eqnarray}
& &\Phi(z)=\sum_{r=1}^4\Phi_r(z)\otimes v_r\  ,\ \ \ \
\Phi^{*}(z)=\sum_{r=1}^4\Phi^{*}_r(z)\otimes v^{*}_r,\\
& &\Psi(z)=\sum_{r=1}^4v_r\otimes\Psi_r(z)\  ,\ \ \ \
\Psi^{*}(z)=\sum_{r=1}^4v^{*}_r\otimes\Psi^{*}_r(z),
\end{eqnarray}
where $v_r\in V_\a$ and $v_r^*\in V_\a^{*S}$. Then we have
\begin{Proposition}: \label{recursion}
The operators $\Phi(z)$ and $\Psi(z)$ with respect to $V_{\a,z}$ are
determined by the components $\Phi_4(z)$ and $\Psi_1(z)$,
respectively. More explicitly,
\bea
&&\Phi_3(z)=-\frac{1}{\sqrt{\a+1}}[\Phi_4(z),f_2]_{q^{-\a-1}},\no\\
&&\Phi_2(z)=[\Phi_3(z),f_1]_q,~~~~
  \Phi_1(z)=-\frac{1}{\sqrt{\a}}[\Phi_2(z),f_2]_{q^{-\a}},\no\\
&&\Psi_2(z)=\frac{1}{\sqrt{\a}}[\Psi_1(z),e_2]_{q^\a},~~~~
  \Psi_3(z)=[\Psi_2(z),e_1]_q,\no\\
&&\Psi_4(z)=\frac{1}{\sqrt{\a+1}}[\Psi_3(z),e_2]_{q^{\a+1}}.
\eea
With respect to $V^{*S}_{\a,z}$, the operators $\Phi^*(z)$ and
$\Psi^*(z)$ are determined by $\Phi^*_1(z)$ and $\Psi^*_4(z)$,
respectively: 
\bea
&&\Phi^*_2(z)=\frac{q^{-\a}}{\sqrt{\a}}[\Phi^*_1(z),f_2]_{q^\a},~~~~
  \Phi^*_3(z)=-q^{-1}[\Phi^*_2(z),f_1]_q,\no\\
&&\Phi^*_4(z)=-\frac{q^{-\a-1}}{\sqrt{\a+1}}[\Phi^*_3(z),f_2]_{q^{\a+1}},\no\\
&&\Psi^*_3(z)=-\frac{q^{\a+1}}{\sqrt{\a+1}}[\Psi^*_4(z),e_2]_{q^{-\a-1}},\no\\
&&\Psi^*_2(z)=-q[\Psi^*_3(z),e_1]_q,~~~~
  \Psi^*_1(z)=\frac{q^\a}{\sqrt{\a}}[\Psi^*_2(z),e_2]_{q^{-\a}}.
\eea
\end{Proposition}

Next we determine the relations between the components $\Phi_4(z),
\Phi^*_1(z), \Psi_1(z), \Psi^*_4(z)$ and the Drinfeld generators. We
have
\begin{Proposition}:
For $\Phi(z)$ associated with $V_{\a,z}$,
\bea
&&[\Phi_4(z), X^{+,i}(w)]=0,~~~i=1,2,\no\\
&&q^{h^i_0}\Phi_4(z)q^{-h^i_0}=q^{-(\a+1)\d_{i2}}\Phi_4(z),\no\\
&&[h^i_n, \Phi_4(z)]=-\d_{i2}q^{(1+\frac{3}{2}k)n}\frac{[(\a+1)n]_q}{n}
  z^n\Phi_4(z),\no\\
&&[h^i_{-n}, \Phi_4(z)]=-\d_{i2}q^{-(1+\frac{1}{2}k)n}\frac{[(\a+1)n]_q}{n}
  z^{-n}\Phi_4(z);\label{phi4}
\eea
for $\Phi^*(z)$ associated with $V^{*S}_{\a,z}$,
\bea
&&[\Phi^*_1(z), X^{+,i}(w)]=0,~~~i=1,2,\no\\
&&q^{h^i_0}\Phi^*_1(z)q^{-h^i_0}=q^{\a\d_{i2}}\Phi^*_1(z),\no\\
&&[h^i_n, \Phi^*_1(z)]=\d_{i2}q^{\frac{3}{2}kn}\frac{[\a n]_q}{n}
  z^n\Phi^*_1(z),\no\\
&&[h^i_{-n}, \Phi^*_1(z)]=\d_{i2}q^{-\frac{1}{2}kn}\frac{[\a n]_q}{n}
  z^{-n}\Phi^*_1(z);\label{phi1}
\eea
for $\Psi(z)$ associated with $V_{\a,z}$,
\bea
&&[\Psi_1(z), X^{-,i}(w)]=0,~~~i=1,2,\no\\
&&q^{h^i_0}\Psi_1(z)q^{-h^i_0}=q^{-\a\d_{i2}}\Psi_1(z),\no\\
&&[h^i_n, \Psi_1(z)]=-\d_{i2}q^{\frac{1}{2}kn}\frac{[\a n]_q}{n}
  z^n\Psi_1(z),\no\\
&&[h^i_{-n}, \Psi_1(z)]=-\d_{i2}q^{-\frac{3}{2}kn}\frac{[\a n]_q}{n}
  z^{-n}\Psi_1(z);\label{psi1}
\eea
and for $\Psi^*(z)$ associated with $V^{*S}_{\a,z}$,
\bea
&&[\Psi^*_4(z), X^{-,i}(w)]=0,~~~i=1,2,\no\\
&&q^{h^i_0}\Psi^*_4(z)q^{-h^i_0}=q^{(\a+1)\d_{i2}}\Psi^*_4(z),\no\\
&&[h^i_n, \Psi^*_4(z)]=\d_{i2}q^{(\frac{1}{2}k-1)n}\frac{[(\a+1)n]_q}{n}
  z^n\Psi^*_4(z),\no\\
&&[h^i_{-n}, \Psi^*_4(z)]=\d_{i2}q^{(-\frac{3}{2}k+1)n}\frac{[(\a+1)n]_q}{n}
  z^{-n}\Psi^*_4(z).\label{psi4}
\eea
\end{Proposition}

\vskip.1in
To obtain bosonized expressions of the intertwining operators,
we introduce the combinations of bosonic oscillators for $m\in{\bf
Z}$,
\bea
&&A^{*}_m=-(a^1_m+\frac{[2m]_q}{[m]_q}a^2_m)q^{\frac{|m|}{2}},\no\\
&&B^{*}_m=-\frac{[\a m]_q}{[(\a+1)m]_q}
  (a^1_m+\frac{[2m]_q}{[m]_q}a^2_m)q^{\frac{|m|}{2}},\no\\
&&\tilde{B}^{*}_m=-(a^1_m+\frac{[2m]_q}{[m]_q}a^2_m)
  q^{-\frac{|m|}{2}}+(b^{13}_m+q^{-|m|}b^{23}_m)
  q^{-\frac{\a}{2}|m|},\no\\
&&\tilde{A}^{*}_m=-\frac{[\a m]_q}{[(\a+1)m]_q}
  (a^1_m+\frac{[2m]_q}{[m]_q}a^2_m)
  q^{-\frac{|m|}{2}}-(q^{|m|}b^{13}_m+b^{23}_m)
  q^{\frac{3\a}{2}|m|},\no\\
&&Q_{A^*}=-Q_{a^1}-2Q_{a^2},~~~~Q_{B^*}=-\frac{\a}{\a+1}
  (Q_{a^1}+2Q_{a^2}),\no\\
&&Q_{\tilde{B}^*}=-Q_{a^1}-2Q_{a^2}+Q_{b^{13}}+Q_{b^{23}},\no\\
&&Q_{\tilde{A}^*}=-\frac{\a}{\a+1}(Q_{a^1}+2Q_{a^2})-Q_{b^{13}}-Q_{b^{23}}.
\eea
For $k=\a$, 
these operators obey the commutation relations, among others,
\bea
&&[A^*_m,h^i_n]=\d_{i2}\d_{m+n,0}\frac{[m]_q[(\a+1)m]_q}{m}
  =[\tilde{A}^*_m, h^i_n],\no\\
&&[B^*_m,h^i_n]=\d_{i2}\d_{m+n,0}\frac{[m]_q[\a m]_q}{m}
  =[\tilde{B}^*_m, h^i_n].
\eea
Then
\begin{Theorem}:
For $k=\a$, the bosonized forms
$\phi_4(z), \phi^*_1(z), \psi_1(z)$ and $\psi^*_4(z)$ of the vertex
operator components $\Phi_4(z), \Phi^*_1(z), \Psi_1(z)$ and $\Psi^*_4(z)$
are given by
\bea
\phi_4(z)&=&:e^{-A^{*}(q^{\a+1}z;-\frac{\a}{2})}:,\no\\
\phi^*_1(z)&=&:e^{B^{*}(q^\a z;-\frac{\a}{2})}:,\no\\
\psi_1(z)&=&:e^{-\tilde{B}^{*}(q^\a z;\frac{\a}{2})}:,\no\\
\psi^*_4(z)&=&:e^{\tilde{A}^{*}(q^{\a-1}z;\frac{\a}{2})}:
    e^{\sqrt{-1}\pi (b^{13}_0+b^{23}_0)}.\label{elementary-v}
\eea
The other components
$\phi_r(z),~ \phi^*_r(z),~\psi_r(z)$ and $\psi^*_r(z)$ are represented
by multiple contour integrals of the Drinfeld currents (c.f. proposition
\ref{recursion}). 
\end{Theorem}

\vskip.1in

Vertex operators (\ref{elementary-v}) 
are referred to as ``elementary $q$-vertex
operators" and are determined solely from their
commutation relations with the bosonized $\Uq$ generators. The
construction is completely independent of which infinite dimensional
modules the vertex operators intertwine. In next section, we shall
clarify on which space these bosonized vertex operators act.

\sect{Fock space and Fock-Wakimoto modules}

In this section we study bosonic Fock space on which
the $\Uq$ generators and the bosonized vertex operators act. 
As we will see, all highest weight modules of $\Uq$ can be embedded in
the bosonic Fock space. Note that $k=\a \neq 0, -1$.

Let $|0>$ be the vacuum
vector, which is defined by $a^i_n |0> = b^{12}_n |0> = b^{13}_n |0> = 
b^{23}_n |0> = c_n |0> = 0$ for $n\geq 0$. Introduce the vector
\bea
&&|\l_{a^1}, \l_{a^2}, \l_{b^{12}}, \l_{b^{13}}, \l_{b^{23}},
  \l_c>\no\\
&&~~~~~~~  =e^{\frac{1}{\a+1}\l_{a^1}Q_{a^1}+\frac{2}{\a+1}
   \l_{a^2}Q_{a^2}+\l_{b^{12}}Q_{b^{12}}+
  \l_{b^{13}}Q_{b^{13}}+\l_{b^{23}}Q_{b^{23}}+\l_cQ_c}|0>,
\eea
which carries the weight $(\frac{\l_{a^1}}{\a+1}, \frac{2\l_{a^2}}{\a+1}, 
  \l_{b^{12}}, \l_{b^{13}},
\l_{b^{23}},\l_c)\in {\bf C}^6$. Denote by
\bea
F_{\frac{1}{\a+1}\l_{a^1},\; \frac{2}{\a+1}
  \l_{a^2},\; \l_{b^{12}},\; \l_{b^{13}},\; \l_{b^{23}}, \;
     \l_c}\no
\eea
the module generated by the creation operators
$a^1_n, a^2_n, b^{12}_n, b^{13}_n, b^{23}_n$ and $c_n ~(n<0)$ over the vector
$|\l_{a^1}, \l_{a^2}, \l_{b^{12}}, \l_{b^{13}}, \l_{b^{23}}, \l_c>$.
Introduce the bosonic Fock space
\bea
&&F_{(\l_{a^1},\l_{a^2},\l_{b^{12}},\l_{b^{13}},\l_{b^{23}},\l_c)}\no\\
&&~~~~~~~~ = \bigoplus_{i,j,l\in{\bf Z}}
 F_{\frac{1}{\a+1}\l_{a^1},\;\frac{2}{\a+1}\l_{a^2},\;\l_{b^{12}}+
   i+j,\; \l_{b^{13}}+j,\; \l_{b^{23}}+l,\; \l_c+i+j}.
\eea
It can be shown that the action of $\Uq$ on this space is closed, i.e.
$\Uq F_* = F_*$ for $*=(\l_{a^1}, \l_{a^2}, \l_{b^{12}}, \l_{b^{13}},
\l_{b^{23}}, \l_c)$. Hence the Fock space $F_*$
constitutes a $\Uq$-module. The elementary $q$-vertex operators are
mapps of the following Fock spaces:
\bea
\phi_r(z),~\psi_r(z)&:& ~~~F_{(\l_{a^1},\,\l_{a^2},\,\l_{b^{12}},\,
   \l_{b^{13}},\,\l_{b^{23}},\,\l_c)}\longrightarrow 
   F_{(\l_{a^1}+\a+1,\,\l_{a^2}+\a+1,\,\l_{b^{12}},\,
   \l_{b^{13}},\,\l_{b^{23}},\,\l_c)},\no\\
\phi^*_r(z),~\psi^*_r(z)&:& ~~~F_{(\l_{a^1},\,\l_{a^2},\,\l_{b^{12}},\,
   \l_{b^{13}},\,\l_{b^{23}},\,\l_c)}\longrightarrow 
   F_{(\l_{a^1}-\a,\,\l_{a^2}-\a,\,\l_{b^{12}},\,
   \l_{b^{13}},\,\l_{b^{23}},\,\l_c)},\label{map-on-F}
\eea
for all $r=1,2,3,4$.

Let us now discuss the emdedding of the highest weight module $V(\l)$ in
the bosonic Fock space $F_*$. We impose the highest weight conditions
on the  vector $|\l_{a^1}, \l_{a^2}, \l_{b^{12}}, \l_{b^{13}},
\l_{b^{23}},\l_c>$,
\bea
&&e_i|\l_{a^1}, \l_{a^2}, \l_{b^{12}}, \l_{b^{13}}, \l_{b^{23}},\l_c>
  =0,\no\\
&&h_i|\l_{a^1}, \l_{a^2}, \l_{b^{12}}, \l_{b^{13}}, \l_{b^{23}},\l_c>
  =\l_i|\l_{a^1}, \l_{a^2}, \l_{b^{12}}, \l_{b^{13}},
  \l_{b^{23}},\l_c> \label{h-condition}
\eea
for all $i=0,1,2$.
Sovling these conditions, we obtain the highest weight vector
$|\b, \g, 0, 0, 0, 0>$, where $\b$ and $\g$
are arbitrary complex parameters. The corresponding highest weight
is $\l_{\b,\,\g}=(\a-\b+2\g)\L_0+
2(\b-\g)\L_1 - \b\L_2$. Thus we have the identification
\beq
|\l_{\b,\g}>=|\b, \g, 0,0,0,0>.
\eeq
Denote by
\bea
F_{(\b,\,\g)}=\bigoplus_{i,j,l\in{\bf Z}}
 F_{\frac{1}{\a+1}\b,\;\frac{2}{\a+1}\g,\;
   i+j,\; j,\; l,\; i+j}
\eea
the Fock space associated to this highest weight vector. It is easy to
see that the $\Uq$ action on the subspace $F_{(\b,\,\g)}$ is
still closed and therefore $F_{(\b,\,\g)}$ is a $\Uq$-module.
Using the highest weight vector $|\l_{\b,\g}>$, 
we construct the level-$\a$ highest
weight module of $\Uq$,
\beq
V(\l_{\b,\g})=\Uq \,|\l_{\b,\g}> .
\eeq
This module is not irreducible in general, but contains a maximal
proper submodule $M(\l_{\b,\g})$ such that $V(\l_{\b,\g}) /
M(\l_{\b,\g})$ yields an irreducible $\Uq$ module. 
It is clear that the module $V(\l_{\b,\g})$ can be embedded in the
bosonic Fock space $F_{(\b,\,\g)}$. Moreover,
from (\ref{map-on-F}) the elementary $q$-vertex
operators are mappings of the Fock spaces:
\bea
\phi_r(z),~\psi_r(z)&:& ~~~F_{(\b,\,\g)}\longrightarrow
   F_{(\b+\a+1,\,\g+\a+1)},\no\\
\phi^*_r(z),~\psi^*_r(z)&:& ~~~F_{(\b,\,\g)}\longrightarrow
   F_{(\b-\a,\,\g-\a)}.\label{map-on-subF}
\eea

However, the Fock space $F_{(\b,\,\g)}$ contains some
redundancies arising from the free bosonic field $c(z;0)$. To see this,
we define the fermionic ghost system $(\eta,\xi)$ of dimension
$(1,0)$,
\beq
\eta(z)=\sum_{n\in{\bf Z}}\eta_n z^{-n-1}=:e^{c(z;0)}:,~~~~~
\xi(z)=\sum_{n\in{\bf Z}}\xi_n z^{-n}=:e^{-c(z;0)}:
\eeq
The mode expansion of $\eta(z)$ and $\xi(z)$ is well defined on
$F_{(\b,\,\g)}$, and the modes satisfy the relations
\beq
\xi_m\xi_n+\xi_n\xi_m=0=\eta_m\eta_n+\eta_n\eta_m,~~~~
  \xi_m\eta_n+\eta_n\xi_m=\d_{m+n,0}.
\eeq
Obviously, $\eta_0\xi_0$ and $\xi_0\eta_0$ qualify as projectors and
so we use them to decompose $F_{(\b,\,\g)}$ into a direct sum of
subspaces
\beq
F_{(\b,\,\g)}=\eta_0\xi_0F_{(\b,\,\g)}
  \oplus \xi_0\eta_0F_{(\b,\,\g)}.
\eeq
$\eta_0\xi_0F_{(\b,\,\g)}$ is referred to as
${\rm Ker}_{\eta_0}$ and $\xi_0\eta_0F_{(\b,\,\g)} = 
F_{(\b,\,\g)}
/ \eta_0\xi_0F_{(\b,\,\g)}$ as 
${\rm Coker}_{\eta_0}$. 
\begin{Proposition}:
$\eta_0$ commutes (or anticommutes) with the action of
$\Uq$. Thus ${\rm Ker}_{\eta_0}$ and
${\rm Coker}_{\eta_0}$ are both $\Uq$-modules.
\end{Proposition}
We are now in a position to consider a restriction of the Fock space
$F_{(\b,\,\g)}$ to a smaller space $\F_{(\b,\,\g)}$, referred to as
the Fock-Wakimoto space.
\begin{Proposition}:
The restricted Fock space
\beq
\F_{(\b,\,\g)}\equiv {\rm Ker}_{\eta_0}F_{(\b,\,\g)}=
  \eta_0\xi_0 F_{(\b,\,\g)}
\eeq
constitutes a Fock-Wakimoto module of $\Uq$.
\end{Proposition}

One can check that $\eta_0 |\l_{\b,\g}>=0$ for any $\b,~\g\in{\bf C}$.
Thus $|\l_{\b,\g}>$ is a $\Uq$ highest weight vector
belonging to the smaller space ${\rm Ker}_{\eta_0}F_{(\b,\,\g)}$. 
It follows that
\begin{Proposition}:
The Fock-Wakimoto module $\F_{(\b,\,\g)}$ is a highest weight
$\Uq$-module with highest weight vector $|\l_{\b,\g}>$ and highest
weight $\l_{\b,\,\g}$.
\end{Proposition}

Using the projection operator $\eta_0\xi_0$, we define the ``projected
$q$-vertex operators" $\tilde{\phi}_r(z),~ \tilde{\phi}^*_r(z),~
\tilde{\psi}_r(z)$ and $\tilde{\psi}^*_r(z)$ as follows
\bea
\tilde{\Theta}(z)=\eta_0\xi_0\Theta(z)\eta_0\xi_0,~~~~
  \Theta(z)=\phi_r(z),~\phi^*_r(z),~\psi_r(z)~{\rm or}~\psi^*_r(z).
  \label{proj-vertex}
\eea
Since $\eta_0$ commutes with the elementary $q$-vertex operators, we
can deduce from (\ref{map-on-subF}) that the projected $q$-vertex
operators are mappings of the highest weight Fock-Wakimoto modules: 
\bea
\tilde{\phi}_r(z)&:& ~~~\F_{(\b,\;\g)}\longrightarrow 
   \F_{(\b+\a+1,\;\g+\a+1)},\no\\
\tilde{\psi}_r(z)&:& ~~~\F_{(\b,\;\g)}\longrightarrow 
   \F_{(\b+\a+1,\;\g+\a+1)},\no\\
\tilde{\phi}^*_r(z)&:& ~~~\F_{(\b,\;\g)}\longrightarrow 
   \F_{(\b-\a,\;\g-\a)},\no\\
\tilde{\psi}^*_r(z)&:& ~~~\F_{(\b,\;\g)}\longrightarrow 
   \F_{(\b-\a,\;\g-\a)}.
\eea

\sect{Screen currents and correlation functions}

Due to the existence of backgound charges, the projected $q$-vertex
operators are not yet the proper bosonizations of the $q$-vertex
operators (\ref{intertwiners}). In this section we construct $q$-screen
currents which  balance the background charges and thus
ensure the nonvanishing of correlation functions of the bosonized
$q$-vertex operators.

Let us introduce the oscillators
\beq
a^{*,i}_m=\frac{[m]_q}{[(k+1)m]_q}a^i_m,~~~~Q_{a^{*,i}}=\frac{1}{k+1}
  Q_{a^i},~~~i=1,2
\eeq
and define the corresponding currents $S^i(z)$ by
\bea
S^i(z)&=&:e^{-a^{*,i}(z;\frac{k+1}{2})}:\tilde{S}^i(z),\\
\tilde{S}^1(z)&=&:e^{-b^{12}(z;0)-b^{12}_-(q^{-1}z)-b^{13}_-(q^{-1}z)
   +b^{23}_-(z)}\,{}_1\partial_z\,e^{-c(q^{-1}z;0)}:\;
   e^{\sqrt{-1}\pi(c_0+b^{12}_0)}\no\\
& &+q\,:e^{b^{13}(z;0)-b^{23}(qz;0)+b^{23}_+(z)}:\;
   e^{-\sqrt{-1}\pi(b^{13}_0+b^{23}_0)},\\
\tilde{S}^2(z)&=&-q^{-1}\,:e^{b^{23}(z;0)}:\;
   e^{-\sqrt{-1}\pi(c_0+b^{12}_0+b^{13}_0+b^{23}_0)}.
\eea
Here we have used the notation 
\beq
{}_k\partial_z f(z)=\frac{f(q^kz)-f(q^{-k}z)}{(q-q^{-1})z}.
\eeq
Then we can verify
\begin{Theorem}:
The currents $S^i(z)$ satisfy the following commutation relations
with the $\Uq$ generators
\bea
&&[h^i_n, S^j(w)]=0,~~~n\in {\bf Z},\no\\
&&[X^{+,i}(z), S^j(w)]=0,\no\\
&&[X^{-,i}(z), S^j(w)]=\d^{ij}\,{}_{k+1}\partial_w\lt(-z^{-1}\,\cdot
  \d(\frac{w}{z}):e^{-a^{*,i}(w;-\frac{k+1}{2})}:\rt).\label{s-current}
\eea
That is, the currents $S^i(z)$ (anti-)commute with the action of $\Uq$ up to
total differences.
The currents $S^i(z)$ are referred to as the q-screen currents of $\Uq$.
\end{Theorem}

For $p\in{\bf C},~|p|<1$ and $s\in{\bf C}-\{0\}$, one defines the
Jackson integral
\beq
\int_0^{s\infty}f(z)d_pz=s(1-p)\sum_{m\in{\bf Z}}f(sp^m)p^m,
\eeq
The Jackson integral enjoys the following property, among others,
\beq
\int_0^{s\infty}f(z)d_pz=\int_0^{s\infty}pf(pz)d_pz,
\eeq
which implies that for $p=q^{2k}$,
\beq
\int_0^{s\infty}{}_k\partial_zf(z)d_pz=0.
\eeq
Note that the right hand side of (\ref{s-current}) is a total
$p=q^{2(k+1)}$ difference. We have
\begin{Corollary}:
The screen charges 
\beq
Q^i=\int_0^{s\infty}S^i(z)d_pz,~~~p=q^{2(k+1)},
\eeq
assuming that the Jackson integrals are convergent, (anti-)commute with all
the generators of $\Uq$.
\end{Corollary}

Since $\eta_0$ commutes with $S^i(z),~i=1,2$, the 
screen charges with $k=\a$ give
rise to the following mappings of the Fock-Wakimoto modules:
\bea
Q^1&:&~~~\F_{(\b,\;\g)}\longrightarrow
    \F_{(\b-1,\;\g)},\\
Q^2&:&~~~\F_{(\b,\;\g)}\longrightarrow
    \F_{(\b,\;\g-\frac{1}{2})}.
\eea

Introduce the screened $q$-vertex operators
\bea
&&\tilde{\phi}_r^{(x_1,\tilde{x}_1)}(z)=(Q^1)^{x_1}(Q^2)^{\tilde{x}_1}
  \tilde{\phi}_r(z),\no\\
&&\tilde{\phi}_r^{*(y_1,\tilde{y}_1)}(z)=(Q^1)^{y_1}(Q^2)^{\tilde{y}_1}
  \tilde{\phi}_r^*(z),\no\\
&&\tilde{\psi}_r^{(x'_1,\tilde{x}'_1)}(z)=(Q^1)^{x'_1}(Q^2)^{\tilde{x}'_1}
  \tilde{\psi}_r(z),\no\\
&&\tilde{\psi}_r^{*(y'_1,\tilde{y}'_1)}(z)=(Q^1)^{y'_1}(Q^2)^{\tilde{y}'_1}
  \tilde{\psi}_r^*(z).
\eea
We are now in a position to state
\begin{Theorem}:
The $q$-vertex operators (\ref{intertwiners}) are bosonized as
\bea
& &\tilde{\Phi}_{\l_{\b,\g}}^{\l_{\b^1_+(x),\g^1_+(\tilde{x})}V}(z)
   =\sum_{r=1}^4\tilde{\phi}^{(x_1,\tilde{x}_1)}_r(z)\otimes v_r ,\no\\
& &\tilde{\Phi}_{\l_{\b,\g}}^{\l_{\b^1_-(y),\g^1_-(\tilde{y})}V^*}(z)
   =\sum_{r=1}^4\tilde{\phi}^{*(y_1,\tilde{y}_1)}_r(z)\otimes v^{*}_r,\no\\
& &\tilde{\Psi}_{\l_{\b,\g}}^{V\,\l_{\b^1_+(x'),\g^1_+(\tilde{x}')}}(z)
   =\sum_{r=1}^4v_r\otimes\tilde{\psi}^{(x'_1,\tilde{x}'_1)}_r(z),\no\\
& &\tilde{\Psi}_{\l_{\b,\g}}^{V^*\,\l_{\b^1_-(y'),\g^1_-(\tilde{y}')}}(z)
   =\sum_{r=1}^4v^*_r\otimes\tilde{\psi}^{*(y'_1,\tilde{y}'_1)}_r(z),
\eea
where
\bea
&&\b^1_+(x)=\b+\a+1-x_1,~~~~
  \g^1_+(\tilde{x})=\g+\a+1-\frac{1}{2}\tilde{x}_1,\no\\
&&\b^1_-(y)=\b-\a-y_1,~~~~
  \g^1_-(\tilde{y})=\g-\a-\frac{1}{2}\tilde{y}_1
\eea
for certain choices of nonnegative integers $x_1,~\tilde{x}_1,~y_1$
and $\tilde{y}_1$. These operators are  intertwiners
of the highest weight $\Uq$-modules:
\bea
\tilde{\Phi}_{\l_{\b,\g}}^{\l_{\b^1_+(x),\g^1_+(\tilde{x})}V}(z)
  &:&~~\F_{(\b,\,\g)}\longrightarrow
  \F_{(\b^1_+(x),\,\g^1_+(\tilde{x}))}\otimes V_{\a,z},\no\\
\tilde{\Phi}_{\l_{\b,\g}}^{\l_{\b^1_-(y),\g^1_-(\tilde{y})}V^*}(z)
  &:&~~\F_{(\b,\,\g)}\longrightarrow
  \F_{(\b^1_-(y),\,\g^1_-(\tilde{y}))}\otimes V_{\a,z}^{*S},\no\\
\tilde{\Psi}_{\l_{\b,\g}}^{V\,\l_{\b^1_+(x'),\g^1_+(\tilde{x}')}}(z)
  &:&~~\F_{(\b,\,\g)}\longrightarrow
  V_{\a,z}\otimes \F_{(\b^1_+(x'),\,\g^1_+(\tilde{x}'))},\no\\
\tilde{\Psi}_{\l_{\b,\g}}^{V^*\,\l_{\b^1_-(y'),\g^1_-(\tilde{y}')}}(z)
  &:&~~\F_{(\b,\,\g)}\longrightarrow
  V_{\a,z}^{*S}\otimes \F_{(\b^1_-(y'),\,\g^1_-(\tilde{y}'))}.
\eea
\end{Theorem}

In the following we compute $N$-point correlation function which is
defined to be the trace of the bosonized $q$-vertex operators over
the $\Uq$-module $\F_{(\b,\,\g)}$, that is
\bea
{\rm Tr}_{\F_{(\b,\,\g)}}\lt(q^{L_0}\Theta_{r_N}(z_N)\cdots
   \Theta_{r_1}(z_1)\rt).
\eea
Here $\Theta_{r_l}(z_l)$ stands for the type I 
$q$-vertex operators 
$\tilde{\phi}_{r_l}^{(x_l,\tilde{x}_l)}(z_l),~~
\tilde{\phi}_{r_l}^{*(y_l,\tilde{y}_l)}(z_l)$
or the type II 
$q$-vertex operators
$\tilde{\psi}_{r_l}^{(x'_l,\tilde{x}'_l)}(z_l),~~
\tilde{\psi}_{r_l}^{*(y'_l,\tilde{y}'_l)}(z_l)$;
$L_0\equiv -d$ is the $q$-Virasoro operator which is bosonized as 
(for $k=\a\neq 0, -1$)
\bea
-L_0&=&\sum_{n>0}\lt(\frac{n^2}{[n]_q[(\a+1)n]_q}\lt(a^1_{-n}a^2_n+
   a^2_{-n}a^1_n+(q^n+q^{-n})a^2_{-n}a^2_n\rt)\rt.\no\\
& &+\lt.\frac{n^2}{[n]_q^2}(b^{12}_{-n}b^{12}_n-b^{13}_{-n}b^{13}_n
  -b^{23}_{-n}b^{23}_n-c_{-n}c_n)\rt)\no\\
& &+\frac{1}{\a+1}\lt(a^1_0a^2_0+(a^2_0)^2+a^1_0+3a^2_0\rt)\no\\
& &+\frac{1}{2}\lt((b^{12}_0)^2-b^{13}_0(b^{13}_0+1)-b^{23}_0
  (b^{23}_0+1)-(c_0)^2\rt).
\eea
The zero mode part of the $a^1_n,\;a^2_n$ oscillators is added to
the $L_0$ operator so that its eigenvalue on $|\l_{\b,\g}>$ is
$\frac{1}{2(\a+1)}(\l_{\b,\g},
\l_{\b,\g}+2\rho)$, where $\rho=\L_0+\L_1+\L_2$.

Let us define the Fock spaces for $s\in {\bf Z}$, 
\beq
F^{(s)}_{(\b,\,\g)}
=\bigoplus_{i,j,l\in {\bf Z}} F
_{\frac{1}{\a+1}\b,\;\frac{2}{\a+1}\g,\;i+j,\;j,\;l,\;i+j+s}.
\eeq
We have  $F^{(0)}_{(\b,\,\g)}=F_{(\b,\,\g)}$. 
It can be shown that $\eta_0,\xi_0$ intertwine various Fock spaces
\begin{eqnarray*}
\eta_0\ \ :\ \ F^{(s)}_{(\b,\,\g)}\longrightarrow
F^{(s+1)}_{(\b,\,\g)},~~~~~
\xi_0\ \ :\ \ F^{(s)}_{(\b,\,\g)}\longrightarrow
F^{(s-1)}_{(\b,\,\g)}.
\end{eqnarray*}
Since $\eta_0^2=0$, we obtain the following BRST complex:
\beq
\cdots~\stackrel{Q_{s-1}=\eta_0}{\longrightarrow}~
    F^{(s)}_{(\b,\,\g)}~\stackrel{Q_{s}=\eta_0}
    {\longrightarrow}~ F^{(s+1)}_{(\b,\,\g)} ~ \stackrel 
    {Q_{s+1}=\eta_0}{\longrightarrow} ~\cdots . 
\eeq
It follows from $\eta_0\xi_0+\xi_0\eta_0=1$, 
that ${\rm Ker}_{Q_{s}}={\rm Im}_{Q_{s-1}}$ for any $s\in{\bf Z}$. 
We have
\begin{Proposition}:
The $N$-point correlation function of the type I vertex operators
\bea
{\rm Tr}_{\F_{(\b,\,\g)}}\lt(q^{L_0}\tilde{\phi}^{(x_N,\tilde{x}_N)}
  _{r_N}(z_N)\cdots\tilde{\phi}^{(x_1,\tilde{x}_1)}_{r_1}(z_1)
  \rt)\neq 0\no
\eea
iff $\a\in{\bf N}$ and $\sum_{i=1}^N x_i=\frac{1}{2}\sum_{i=1}^N
\tilde{x}_i=N(\a+1)$. For such $\a$ and $x_i,~ \tilde{x}_i$, 
the above trace is given by 
\beq
\sum_{s=1}^\infty (-1)^{s+1}{\rm Tr}_{F^{(-s)}_{(\b,\,\g)}}\lt(
q^{L_0}(Q^1)^{x_N}(Q^2)^{\tilde{x}_N}\phi_{r_N}(z_N)\cdots
(Q^1)^{x_1}(Q^2)^{\tilde{x}_1}\phi_{r_1}(z_1)\rt).
\eeq
Similarly, the $N$-point correlator of the type II vertex operators
\bea
&&{\rm Tr}_{\F_{(\b,\,\g)}}\lt(q^{L_0}\tilde{\psi}^{(x'_N,\tilde{x}'_N)}
  _{r_N}(z_N)\cdots\tilde{\psi}^{(x'_1,\tilde{x}'_1)}_{r_1}(z_1)\rt)\no\\
&&~~~~~~~=\sum_{s=1}^\infty (-1)^{s+1}{\rm Tr}_{F^{(-s)}_{(\b,\,\g)}}\lt(
q^{L_0}(Q^1)^{x'_N}(Q^2)^{\tilde{x}'_N}\psi_{r_N}(z_N)\cdots
(Q^1)^{x'_1}(Q^2)^{\tilde{x}'_1}\psi_{r_1}(z_1)\rt)\no\\
\eea
is nonvanishing
iff $\a\in{\bf N}$ and $\sum_{i=1}^N x'_i=\frac{1}{2}\sum_{i=1}^N
\tilde{x}'_i=N(\a+1)$. 
\end{Proposition}

We now consider the $N$-point correlation function involing also dual
vertex operators,
\beq
{\rm Tr}_{\F_{(\b,\,\g)}}\lt(q^{L_0}\tilde{\phi}^{*(y_N,\tilde{y}_N)}
  _{r_N}(z_N)\cdots\tilde{\phi}^{*(y_{l+1},\tilde{y}_{l+1})} _{r_{l+1}}
  (z_{l+1})\tilde{\phi}^{(x_l,\tilde{x}_l)}_{r_l}(z_l)\cdots
  \tilde{\phi}^{(x_1,\tilde{x}_1)}_{r_1}(z_1)\rt).
  \label{n-point-dual}
\eeq
Then we have
\begin{Proposition}:
For $\a\in{\bf N}$, (\ref{n-point-dual}) is non-zero iff
$\sum_{i=1}^lx_i+\sum_{i=l+1}^Ny_i=\frac{1}{2}(\sum_{i=1}^l
  \tilde{x}_i+\sum_{i=l+1}^N\tilde{y}_i)=(2l-N)\a+l$. And for 
$\a\not\in{\bf N}$ it is nonvanishing iff $N$ is even, i.e. $N=2L$,
and $l=L=\sum_{i=1}^Lx_i+\sum_{i=L+1}^Ny_i=\frac{1}{2}(\sum_{i=1}^L
\tilde{x}_i+\sum_{i=L+1}^N\tilde{y}_i)$. In both case, the trace
(\ref{n-point-dual}) can be written as the following unified formula
\bea
(\ref{n-point-dual})&=&\sum_{s=1}^\infty (-1)^{s+1}{\rm Tr}
  _{F^{(-s)}_{(\b,\,\g)}}\lt(q^{L_0}(Q^1)^{y_N}(Q^2)^{\tilde{y}_N}
  \phi^*_{r_N}(z_N)\cdots(Q^1)^{y_{l+1}}(Q^2)^{\tilde{y}_{l+1}}
  \phi^*_{r_{l+1}}(z_{l+1})\rt.\no\\
& &~~~~~~\lt.(Q^1)^{x_{l}}(Q^2)^{\tilde{x}_{l}}\phi_{r_{l}}(z_{l})
  \cdots (Q^1)^{x_{1}}(Q^2)^{\tilde{x}_{1}}\phi_{r_{1}}(z_{1})\rt).
\eea
An integral formula for the $N$-point functions of type II (dual) vertex 
operators can be written down in a similar way, which we omit.
\end{Proposition}

\vskip.3in
\noindent {\large\bf Acknowledgements}
\vskip.1in
This work has been financially supported by the Australian Research
Council large, small and QEII fellowship grants. Discussions with
Wen-Li Yang are gratefully acknowledged.


\end{document}